\documentclass[11pt]{amsart}
\usepackage{amsmath,amssymb,amscd,amsfonts,verbatim}
\usepackage[mathscr]{eucal}
\newtheorem{theorem}{Theorem}
\newtheorem{corollary}{Corollary}
\newtheorem{lemma}{Lemma}
\theoremstyle{definition}
\newtheorem{definition}{Definition}

\newcommand{\IR}{\mathbb R}
\newcommand{\IQ}{\mathbb Q}
\newcommand{\II}{\mathbb I}

\newcommand{\w}{\omega}
\newcommand{\HD}{\mathit{HD}}
\newcommand{\D}{\mathsf{D}}
\newcommand{\la}{\langle}
\newcommand{\ra}{\rangle}

\newcommand{\pr}{\mathrm{pr}}
\newcommand{\e}{\varepsilon}
\newcommand{\diam}{\mathrm{diam}}
\newcommand{\HH}{\mathcal H}
\newcommand{\add}{\mathrm{add}}
\newcommand{\cl}{\mathrm{cl}}
\newcommand{\A}{{\mathcal A}}

\newcommand{\mc}{\mathscr }
\newcommand{\Class}{\sigma\mathfrak C}

\newcommand{\ind}{\mathrm{ind}}

\begin{document}

\title[The topology of systems of hyperpsaces]{The topology of systems of hyperspaces determined by dimension functions}
\author{Taras Banakh}
\address{Department of Mathematics, Ivan Franko University in Lviv, Ukraine, and
Instytut Matematyki, Jan Kochanowski University in Kielce, Poland}
\email{tbanakh@yahoo.com}
\author{Natalia Mazurenko}
\address{Precarpatian National University, Ivano-Frankivsk, Ukraine}
\email{mnatali@ukr.net}
\subjclass{57N17; 54B20; 54F45; 55M10; 54F65; 28A80}
\keywords{Hyperspace, dimension function, Hausdorff dimension, Hilbert cube, absorbing system.}
\thanks{This work was supported by State Fund of fundamental research of Ukraine, project
F25.1/099.}

\begin{abstract}Given a non-degenerate Peano continuum $X$, a dimension function $\D:2^X_*\to[0,\infty]$ defined on the family $2^X_*$ of compact subsets of $X$, and a subset $\Gamma\subset[0,\infty)$, we recognize the topological structure of the system $\la 2^X,\D_{\le\gamma}(X)\ra_{\alpha\in\Gamma}$, where $2^X$ is the hyperspace of non-empty compact subsets of $X$ and
$\D_{\le\gamma}(X)$ is the subspace of $2^X$, consisting of non-empty compact subsets $K\subset X$ with $\D(K)\le\gamma$.
\end{abstract}
\maketitle

\section{Introduction}

The problem of topological characterization (identification) of
topological objects is a central problem in topology.
A classical result of this sort is the Curtis-Schori
Theorem \cite{CS} asserting that for each non-degenerate Peano continuum
$X$ the hyperspace $2^X$ of non-empty compact subsets of $X$ endowed with the Vietoris topology is homeomorphic to the Hilbert cube $Q=[-1,1]^\w$. At bit later, D.Curtis \cite{curtis} characterized topological spaces $X$ whose hyperspace $2^X$ is homeomorphic to the pseudointerior $s=(-1,1)^\w$ of the Hilbert cube as connected locally connected Polish nowhere locally compact spaces.

In \cite{DR} T.Dobrowolski and L.Rubin recognized the topology of the subspace $\dim_{\le n}(Q)\subset 2^Q$ consisting of compact subsets of $Q$ having covering dimension $\le n$. They constructed a homeomorphism $h:2^Q\to Q^\w$ such that $h(\dim_{\le n}(Q))=Q^n\times s^{\w\setminus n}$ for all $n=\{0,\dots,n-1\}\in \w$. In this case it is said that the system $\la 2^Q,\dim_{\le n}(Q)\ra_{n\in\w}$ is homeomorphic to the system $\la Q^\w,Q^n\times s^{\w\setminus n}\ra_{n\in\w}$.

This result was later generalized by H.Gladdines \cite{glad} to products of Peano continua. Finally,  R.Cauty \cite{coty} has characterized spaces $X$ for which the system $\la 2^X,\dim_{\le n}(X)\ra_{n\in\w}$ is homeomorphic to $\la Q^\w,Q^n\times s^{\w\setminus n}\ra_{n\in\w}$ as Peano continua whose any non-empty open subset contains compact subsets of arbitrary high finite dimension.

In \cite{maz} given a metric space $X$ the second author initiated the study of the subspace $\HD_{\le
\gamma}(X)\subset 2^X$ of compact subsets of $X$ whose Hausdorff dimension is $\le \gamma$. Unlike the
(integer-valued) topological dimension, the Hausdorff dimension of a metric compactum can take on any
non-negative real value $\gamma$. So, the system $\la 2^X,\HD_{\le\gamma}(X)\ra_{\gamma\in[0,\infty)}$ that
naturally appears in this situation is uncountable. In \cite{MazMS} it was proved that for a finite-dimensional
cube $X=[0,1]^n$ the system $\la 2^X,\HD_{\le \gamma}(X)\ra_{\gamma\in[0,n)}$ is homeomorphic to the system $\la
Q^\IQ,Q^{\IQ_{\le\gamma}}\times s^{\IQ_{>\gamma}}\ra_{\gamma\in[0,n)}$ (by $\mathbb Q$ we denote the space of
rational numbers). Here for a subset $A\subset\IR$ and a real number $\gamma$ we put
$$
\begin{gathered}A_{\le \gamma}=\{a\in A:a\le\gamma\},\quad A_{\ge\gamma}=\{a\in A:a\ge\gamma\}\\
A_{<\gamma}=\{a\in A:a<\gamma\},\quad A_{>\gamma}=\{a\in A:a>\gamma\}
\end{gathered}
.$$

Both the (topological) covering dimension and the (metric) Hausdorff dimension are particular cases of dimension functions defined as follows.

\begin{definition} A function $\D: 2_*^X\to[0,\infty]$ defined on the family $2^X_*$
of compact subsets of a topological space $X$ is called a {\em dimension function} if:
\begin{enumerate}
\item $\D(\emptyset)=0$;
\item $\D$ is {\em 
monotone} in the sense that $\D(A)\le\D(B)$ for any compact subsets $A\subset B$ of $X$; \item $\D$ is {\em finitely additive} in
the sense that $\D(F\cup A\cup B)\le\max\{\D(A),\D(B)\}$ for any finite subset $F\subset X$ and disjoint compact subsets $A,B\subset X$; 
\item $\D$ is {\em
$\w$-additive} in the sense that each non-empty open subset $U\subset X$ contains non-empty open sets
$U_n\subset U$, $n\in\w$, such that each compact subset $K\subset \cl_X(\bigcup_{n\in\w}U_n)$ has dimension
$\D(K)\le\sup_{n\in\w}\D(K\cap \overline{U}_n)$.
\end{enumerate}
\end{definition}

\smallskip


Given a dimension function $\D:2_*^X\to[0,\infty]$ on $X$ and a subset $\Gamma\subset[0,\infty)$,
for every $\gamma\in\Gamma$ consider  the subspace $$\D_{\le\gamma}(X)=\{F\in 2^X:\D(F)\le\gamma\}$$ in the
hyperspace $2^X$. Our aim is to recognize the topological structure of the system $\la
2^X,D_{\le\gamma}(X)\ra_{\gamma\in\Gamma}$.

In the sequel, by a $\Gamma$-system $\la X,X_\gamma\ra_{\gamma\in\Gamma}$ we shall  understand a pair consisting
of a set $X$ and a family $\la X_\gamma\ra_{\gamma\in\Gamma}$ of subsets of $X$,  indexed by the elements of an
index set $\Gamma$. Two $\Gamma$-systems $\la X,X_\gamma\ra_{\gamma\in\Gamma}$ and $\la
Y,Y_\gamma\ra_{\gamma\in\Gamma}$ are {\em homeomorphic} if there is a homeomorphism $h:X\to Y$ such that
$h(X_\gamma)=Y_\gamma$ for all $\gamma\in\Gamma$.


The following theorem describes the topological structure  of the $\Gamma$-system $\la
2^X,\D_{\le\gamma}(X)\ra_{\gamma\in\Gamma}$ for a dimension function $\D:2^X_*\to[0,\infty]$ taking values
in the half-line with attached infinity  (that is assumed to be larger than any real number). In that theorem we
shall refer to the subsets $(\gamma]_\Gamma$  defined for $\Gamma\subset\IR$ and $\gamma\in\Gamma$ as follows:
$$(\gamma]_\Gamma=\begin{cases}
(\gamma,\inf(\Gamma_{>\gamma})]&\mbox{if \ $\gamma<\inf(\Gamma_{>\gamma})$;}\\
(\sup(\Gamma_{<\gamma}),\gamma]&\mbox{if \ $\Gamma\ni\sup(\Gamma_{<\gamma})<\gamma=\inf(\Gamma_{>\gamma})$;}\\
[\sup(\Gamma_{<\gamma}),\gamma]&\mbox{in all other cases.}\\
\end{cases}
$$
In this definition we assume that $\sup(\emptyset)=-\infty$ and $\inf(\emptyset)=+\infty$.


\begin{theorem}\label{main} Let $X$ be a topological space and $\D:2_*^X\to[0,\infty]$ be a dimension function. For every subset $\Gamma\subset [0,\infty)$ the $\Gamma$-system $\la 2^X,\D_{\le \gamma}(X)\ra_{\gamma\in\Gamma}$ is homeomorphic to the $\Gamma$-system $\la Q^\IQ,Q^{\IQ_{\le\gamma}}\times s^{\IQ_{>\gamma}}\ra_{\gamma\in\Gamma}$ if and only if
\begin{enumerate}
\item $X$ is a non-degenerate Peano continuum,
\item each subspace $\D_{\le \gamma}(X)$, $\gamma\in\Gamma$, is of type $G_\delta$ in $2^X$, and
\item each non-empty open set $U\subset X$ for every $\gamma\in\Gamma$ contains a compact subset $K\subset U$ with $\D(K)\in(\gamma]_\Gamma$.
\end{enumerate}
\end{theorem}

First, we apply this theorem to integer-valued dimension functions. We identify each natural number
$n$ with the set $\{0,\dots,n-1\}$. Also we put $\overline{\w}=\w\cup\{\w\}$.

\begin{corollary} Let $X$ be a topological space and $\D:2_*^X\to\overline{\w}$ be a dimension function.
For every $n\in\overline{\w}$ the $n$-system $\la 2^X,\D_{\le k}(X)\ra_{k\in n}$ is homeomorphic to the $n$-system
$\la Q^\w,Q^k\times s^{\w\setminus k}\ra_{k\in n}$ if and only if
\begin{enumerate}
\item $X$ is a non-degenerate Peano continuum, \item each subspace $\D_{\le k}(X)$, $k\in n$, is
of type $G_\delta$ in $2^X$, and \item each non-empty open set $U\subset X$ for every $k\in n$ contains a
compact subset $K\subset U$ with $\D(K)=k$.
\end{enumerate}
\end{corollary}

The covering dimension $\dim$ and the cohomological dimension $\dim_G$ for an arbitrary Abelian group $G$ are examples of integer-valued dimension functions. Therefore Corollary 1 implies the following theorem of R. Cauty \cite{coty} that was mentioned  above.

\begin{theorem}[Cauty]\label{cautyteo}
For any non-degenerate Peano continuum $X$ the $\omega$-systems $\la 2^X,\dim_{\le n}(X)\ra_{n\in\omega}$ is homeomorphic to $\la Q^\w,Q^n\times s^{\w\setminus
n}\ra_{n\in\omega}$ if and only if each non-empty open set $U\subset X$ contains an compact subset of arbitrary finite dimension.
\end{theorem}

In \cite{coty} R. Cauty notices, that this theorem holds also for the cohomological dimension $\dim_G$ or any other dimension function in the sense of \cite{DR}. It does not demand any modifications of arguments in the proof.

\bigskip

Applying  Theorem~\ref{main} to the half-interval $\Gamma=[0,b)\subset[0,\infty)$, we obtain:

\begin{corollary}\label{cor2} Let $X$ be a topological space and $\D:2_*^X\to[0,\infty]$ be a dimension function.
For every $b\in[0,\infty]$ the $[0,b)$-system $\la 2^X,\D_{\le \gamma}(X)\ra_{\gamma\in[0,b)}$ is homeomorphic
to the $[0,b)$-system $\la Q^\IQ,Q^{\IQ_{\le\gamma}}\times s^{\IQ_{>\gamma}}\ra_{\gamma\in[0,b)}$ if and only if
\begin{enumerate}
\item $X$ is a non-degenerate Peano continuum, 
\item each subspace $\D_{\le \gamma}(X)$, $\gamma\in[0,b)$, is of
type $G_\delta$ in $2^X$, and \item each non-empty open set $U\subset X$ for  every $\gamma\in[0,b)$ contains a
compact subset $K\subset U$ with $\D(K)=\gamma$.
\end{enumerate}
\end{corollary}

Applying Corollary~\ref{cor2} to the Hausdorff dimension $\dim_H$ we obtain the following theorem whose partial case for $X=\II^n$ was proved in \cite{MazMS}.

\begin{theorem}\label{haus} For a number $b\in(0,\infty]$ and a non-degenerate metric Peano continuum $X$ the system
$\la 2^X,\HD_{\le\gamma}(X)\ra_{\gamma\in[0,b)}$ is homeomorphic to the system $\la
Q^\IQ,Q^{\IQ_{\le\gamma}}\times s^{\IQ_{>\gamma}}\ra_{\gamma\in[0,b)}$  if and only if each non-empty open
subset $U\subset X$ has  Hausdorff dimension $\dim_H(U)\ge b$.
\end{theorem}

To derive this theorem from Corollary~\ref{cor2}, we need to check the conditions (2) and (3) for the Hausdorff dimension. The condition (2) was establised in \cite{maz} while (3) follows from the subsequent Mean Value Theorem for Hausdorff dimension, which will be proved in Section~\ref{s6}.

\begin{theorem}\label{dimH} Let $X$ be a separable complete metric space $X$. For every non-negative real number $d<\dim_H(X)$ the space $X$ contains a compact subset $K\subset X$ of Hausdorff dimension $\dim_H(K)=d$.
\end{theorem}

A similar Mean Value Theorem holds for topological dimension: each regular space $X$ with finite inductive dimension $\ind(X)$ contains a closed subspaces of any dimension $k\le\ind(X)$, see \cite[1.5.1]{En}. However, (in contrast to the Hausdorff dimension)  this theorem does not hold for infinite-dimensional spaces: there is an infinite-dimensional compact metrizable space $X$ containing no subspace of positive finite dimension \cite[5.2.23]{En}.

\section{Absorbing systems in the Hilbert cube}\label{s1}

Theorem~\ref{main} is proved by the technique of absorbing systems created and developed in \cite{dij}, \cite{glad}. So, in this section we start by recalling some basic information related to absorbing systems.

From now on all topological spaces  are metrizable and separable, all maps are continuous. By $\II$ we denote
the unit interval $[0,1]$, by $\mathbb Q$ the space of rational numbers, by $Q=[-1,1]^\w$ the {\em Hilbert
cube}, by $s=(-1,1)^\w$ its {\em pseudointerior} and by $B(Q)$ its {\em pseudoboundary}. By {\em a Hilbert cube}
we understand any topological space homeomorphic to the Hilbert cube $Q$. In particular, for each at most countable
set $A$ the power $Q^A$ is a Hilbert cube; $B(Q^A)=Q^A\setminus s^A$ will stand for its pseudoboundary.

Given two maps $f,g:X\to Y$ and a cover ${\mathcal U}$ of $Y$ we write $(f,g)\prec{\mathcal U}$ and say that $f,g$ are {\em ${\mathcal U}$-near} if for every point $x\in X$ there is a set $U\in{\mathcal U}$ such that $\{f(x),g(x)\}\subset U$.

A closed subset $A$ of an ANR-space $X$ is a called a {\em $Z$-set} if for each map $f:Q\to X$ and an open cover ${\mathcal U}$ of $X$ there is a map $g:Q\to X\setminus A$ such that $(f,g)\prec{\mathcal U}$.  A subset $A\subset X$ is called a {\em $\sigma Z$-set} if $A$ can be written as the countable union of $Z$-sets. It is known \cite{jan} that a closed $\sigma Z$-set in a Polish ANR-space is a $Z$-set. An embedding $f:K\to X$ is called a {\em $Z$-embedding} if the image $f(K)$ is a $Z$ set in $X$.

It is well-known that each map $f:K\to Q$ defined on a compact space can be approximated by $Z$-embeddings, see \cite{Chap}, \cite{jan}.
\smallskip

Let $\Gamma$ be a set. By a {\em $\Gamma$-system} $\mc X=\la X,X_\gamma\ra_{\gamma\in\Gamma}$ we shall understand a pair consisting of a space $X$ and an indexed
collection $\la X_\gamma\ra_{\gamma \in \Gamma}$ of subsets of $X$. Given a map $f:Z\to X$ and a set  $K\subset X$ let $f^{-1}(\mc X)=\la f^{-1}(X),f^{-1}(X_\gamma)\ra_{\gamma\in\Gamma}$ and $K\cap\mc X=\la K\cap X,K\cap X_\gamma\ra_{\gamma\in\Gamma}$.

From now on, ${\mathfrak C}_\Gamma$ is a fixed class of $\Gamma$-systems.

 Generalizing the standard concept of a strongly universal pair \cite[\S1.7]{BRZ} to $\Gamma$-systems we get an important notion of a strongly ${\mathfrak C}_\Gamma$-universal $\Gamma$-system.

\begin{definition} A $\Gamma$-system $\mc X=\la X,X_\gamma\ra_{\gamma\in\Gamma}$ is defined to {\em strongly $\mathfrak C_\Gamma$-universal} if for any open cover ${\mathcal U}$ on $X$, any $\Gamma$-system $\mc C=\la C,C_\gamma\ra_{\gamma\in\Gamma}\in{\mathfrak C}_\Gamma$, and a map $f:C\to X$ whose restriction $f|B:B\to X$ to a closed subset $B\subset C$ is a $Z$-embedding with $ (f|B)^{-1}(\mc X)=B\cap\mc C$ there exists a $Z$-embedding $\tilde f:C\to X$ such that $(\tilde f,f)\prec{\mathcal U}$, $\tilde f|B=f|B$, and $\tilde f^{-1}(\mc X)=\mc C$.
\end{definition}

The strong universality is the principal ingredient in the notion of a ${\mathfrak C}_\Gamma$-absorbing system, generalizing the notion of an absorbing pair, see \cite[\S1.6]{BRZ}.

\begin{definition}\label{absorb}
A $\Gamma$-system $\mc X=\la X,X_\gamma\ra_{\gamma\in\Gamma}$ is defined to {\em  $\mathfrak C_\Gamma$-absorbing} if
\begin{enumerate}
\item[i)] $\mc X$ is strongly ${\mathfrak C}_\Gamma$-universal;
\item[ii)] there is a sequence $\la Z_n\ra_{n\in\w}$ of $Z$-sets in $X$ such that
$\bigcup_{\gamma\in\Gamma}X_\gamma\subset\bigcup_{n\in\w}Z_n$ and  $Z_n\cap \mc X\in{\mathfrak C}_\Gamma$ for all $n\in\w$.
\end{enumerate}
\end{definition}

A remarkable feature of ${\mathfrak C}_\Gamma$-absorbing system in the Hilbert cube is their topological equivalence. We define two $\Gamma$-systems $\la X,X_\gamma\ra_{\gamma\in\Gamma}$ and $\la Y,Y_\gamma\ra_{\gamma\in\Gamma}$ to be {\em homeomorphic} if there is a homeomorphism $h:X\to Y$ such that $h(X_\gamma)=Y_\gamma$ for $\gamma\in\Gamma$.

The following Uniqueness Theorem can be proved by analogy with Theorem 1.7.6 from \cite{BRZ}.

\begin{theorem}\label{unique} Two ${\mathfrak C}_\Gamma$-absorbing $\Gamma$-systems $\la X,X_\gamma\ra_{\gamma\in\Gamma}$ and $\la Y,Y_\gamma\ra_{\gamma\in\Gamma}$ are homeomorphic provided $X$ and $Y$ are homeomorphic to a manifold modeled on $Q$ or $s$.
\end{theorem}

By a {\em manifold} modeled on a space $E$ we understand a metrizable separable space $M$ whose any point has an
open neighborhood homeomorphic to an open subset of the model space $E$.

\section{Characterizing model absorbing systems}

In this section, given a subset $\Gamma\subset\IR$ we characterize the topology of the model $\Gamma$-system
$\la Q^\IQ,Q^{\IQ_{\le\gamma}}\times s^{\IQ_{>\gamma}}\ra_{\gamma\in\Gamma}$. In fact, it will be more
convenient to work with the complementary $\Gamma$-system
$$\Sigma_\Gamma=\la Q^\IQ,Q^{\IQ_{\le\gamma}}\times B(Q^{\IQ_{>\gamma}})\ra_{\gamma\in\Gamma},$$
where $B(Q^{\IQ_{>\gamma}})=Q^{\IQ_{>\gamma}}\setminus s^{\IQ_{>\gamma}}$.
We shall prove that the latter system is $\Class_\Gamma$-absorbing for a suitable class $\Class_\Gamma$ of $\Gamma$-systems.

Let $\Gamma\subset\IR$.
Let us define a $\Gamma$-system $\la A,A_\gamma\ra_{\gamma\in\Gamma}$ to be
\begin{itemize}
\item {\em $\sigma$-compact} if the space $A$ is compact while all subspaces $A_\gamma$, $\gamma\in\Gamma$, are $\sigma$-compact;
\item {\em $\inf$-continuous} if $A_\gamma=\bigcup_{\beta\in B}A_\beta$ for any subset $B\subset\Gamma$ with $\inf B=\gamma\in\Gamma$.
\end{itemize}

By $\Class_\Gamma$ we shall denote the class of $\sigma$-compact $\inf$-continuous $\Gamma$-systems. Let us observe that each $\Gamma$-system $\la A,A_\gamma\ra_{\gamma\in\Gamma}\in\Class_\Gamma$ is decreasing. Indeed, for any real numbers $\alpha<\beta$ in $\Gamma$ the equality $\alpha=\inf\{\alpha,\beta\}$ implies $A_\alpha=A_\alpha\cup A_\beta\supset A_\beta$.

Each $\Gamma$-system $\A=\la A,A_\gamma\ra_{\gamma\in\Gamma}\in\Class_\Gamma$ can be extended to the $\IR$-system $\tilde\A=\la A,\tilde A_\gamma\ra_{\gamma\in\IR}\in\Class_\IR$ consisting of the sets
$$\tilde A_\gamma=\begin{cases}
\bigcup_{\alpha\in\Gamma_{\ge\gamma}}A_\alpha&\mbox{ if \  $\sup(\Gamma_{<\gamma})\notin\Gamma$ or $\gamma=\inf(\Gamma_{\ge\gamma})$ };\\
A_\alpha&\mbox{ if \ $\alpha=\sup(\Gamma_{<\gamma})\in\Gamma$ and $\gamma<\inf(\Gamma_{\ge\gamma})$ },
\end{cases}
$$
indexed by real numbers $\gamma$.

\begin{lemma}\label{enlarge} The $\IR$-system $\tilde\A=\la A,\tilde A_\gamma\ra_{\gamma\in\IR}$ is $\sigma$-compact, $\inf$-continuous and extends the $\Gamma$-system $\A=\la A,A_\gamma\ra_{\gamma\in\Gamma}$ in the sense that $\tilde A_\gamma=A_\gamma$ for all $\gamma\in\Gamma$.
\end{lemma}

\begin{proof} To see that the $\IR$-system $\tilde \A$ is $\sigma$-compact, fix any real number $\gamma$.
The set $\tilde A_\gamma$ is clearly $\sigma$-compact if $\tilde A_\gamma=A_\alpha$ for some $\alpha\in\Gamma$.
So, we assume that $\tilde A_\gamma\ne A_\alpha$ for all $\alpha\in\Gamma$.
In this case $\tilde A_\gamma=\bigcup_{\alpha\in \Gamma_{\ge\gamma}}A_\alpha$ and $\inf\Gamma_{\ge\gamma}\notin\Gamma$. Choose any countable dense subset $D\subset\Gamma$ and observe that $\inf D_{\ge\gamma}=\inf \Gamma_{\ge\gamma}$ and hence $$\tilde A_\gamma=\bigcup_{\alpha\in \Gamma_{\ge\gamma}}A_\alpha=\bigcup_{\alpha\in D_{\ge\gamma}}A_\alpha$$ is $\sigma$-compact, being the countable union of $\sigma$-compact spaces $A_\alpha$, $\alpha\in D_{\ge \gamma}$.
\smallskip

Observe that for every $\gamma\in\Gamma$ we get
$\gamma=\inf(\Gamma_{\ge\gamma})$ and hence
$A_\gamma\subset \bigcup_{\alpha\in\Gamma_{\ge\gamma}}A_\alpha=\tilde A_\gamma$.
The reverse inclusion $\tilde A_\gamma=\bigcup_{\alpha\in\Gamma_{\ge\gamma}}A_\alpha\subset A_\gamma$ follows from the decreasing property of the $\gamma$-system $\mathcal A$. Thus $A_\gamma=\tilde A_\gamma$, which means that the $\IR$-system $\tilde\A$ extends the $\Gamma$-system $\mathcal A$.
\smallskip

Next, we prove that the $\IR$-system $\tilde \A$ is decreasing. Given two real numbers $\beta<\gamma$, we need to show that $\tilde A_\beta\supset \tilde A_\gamma$. We consider four cases:
\smallskip

1) Both $\beta$ and $\gamma$ satisfy the first case of the definition of $\tilde A_\beta$ and $\tilde A_\gamma$:$$\big(\sup(\Gamma_{<\beta})\notin\Gamma\mbox{ or }\beta=\inf(\Gamma_{\ge\beta})\big)\mbox{ and } \big(\sup(\Gamma_{<\gamma})\notin\Gamma\mbox{ or }\gamma=\inf(\Gamma_{\ge\gamma})\big).$$ In this case $\beta<\gamma$   implies $\Gamma_{\ge\beta}\supset\Gamma_{\ge\gamma}$ and thus $$\tilde A_\beta=\bigcup_{\alpha\in\Gamma_{\ge\beta}}A_\alpha\supset \bigcup_{\alpha\in\Gamma_{\ge\gamma}}A_\gamma=\tilde A_\gamma.$$

 2) The element $\beta$ satisfies the first case of the definition of $\tilde A_\beta$ while $\gamma$ satisfies the second case: $$\big(\sup(\Gamma_{<\beta})\notin\Gamma\mbox{  or }\beta=\inf(\Gamma_{\ge\beta})\big)\mbox{ and } \alpha=\sup(\Gamma_{<\gamma})\in\Gamma\mbox{  and }\gamma<\inf(\Gamma_{\ge\gamma}).$$

In this case $\beta\le \alpha$. Indeed, assuming conversely that
$\alpha<\beta$, we get $\Gamma_{<\beta}=\Gamma_{<\gamma}$ and thus $\alpha=\sup(\Gamma_{<\beta})\in\Gamma$, which implies that $\beta=\inf(\Gamma_{\ge\beta})$. In this case, $\alpha=\sup(\Gamma_{<\gamma})\ge\beta$, which is a contradiction. So, $\beta\le\alpha$ and then $\alpha\in\Gamma_{\ge\beta}$ and $\tilde A_\beta\supset A_\alpha=\tilde A_\gamma$.

3) The element $\beta$ satisfies the second case of the definition of $\tilde A_\beta$ while $\gamma$ satisfies
the first one: $$\alpha=\sup(\Gamma_{<\beta})\in\Gamma\mbox{  and }\beta<\inf(\Gamma_{\ge\beta})\mbox{  and }
\big(\sup(\Gamma_{<\gamma})\notin\Gamma\mbox{  or } \gamma=\inf(\Gamma_{\ge\gamma})\big).$$ In this case $$\tilde
A_\beta=A_\alpha\supset\bigcup_{\delta\in\Gamma_{\ge\gamma}}A_\delta=\tilde A_\gamma.$$

4) Both $\beta$ and $\gamma$ satisfy the second case of the definition of $\tilde A_\beta$ and $\tilde A_\gamma$:
$$\alpha_\beta=\sup(\Gamma_{<\beta})\in\Gamma,\;\beta<\inf(\Gamma_{\ge\beta}),\;\alpha_\gamma=\sup(\Gamma_{<\gamma})\in\Gamma,\;\gamma<\inf(\Gamma_{\ge\gamma}).$$ In this case $\alpha_\beta\le\alpha_\gamma$ and
$\tilde A_\beta=A_{\alpha_{\beta}}\supset A_{\alpha_\gamma}=\tilde A_\gamma$. This completes the proof of the
decreasing property of the $\IR$-system $\tilde A$.
\smallskip

Finally, we show that the $\IR$-system $\tilde \A$ is $\inf$-continuous. Fix any real number $\gamma$ and a subset $B\subset \IR$ with $\gamma=\inf B$. We need to check that $\tilde A_\gamma=\bigcup_{\beta \in B}\tilde A_\beta$.
The decreasing property of $\tilde A$ guarantees that $\tilde A_\gamma\supset \bigcup_{\beta \in B}\tilde A_\beta$. It remains to prove the reverse inclusion, which is trivial if $\gamma\in B$. So, we assume that $\gamma\notin B$.
Two cases are possible:
\smallskip

1. $\sup(\Gamma_{<\gamma})\notin \Gamma$ or $\gamma=\inf(\Gamma_{\ge\gamma})$. In this case $\tilde A_\gamma=\bigcup_{\alpha\in\Gamma_{\ge\gamma}}A_\alpha$.  We consider three subcases:
\smallskip

1a) If $\gamma=\inf(\Gamma_{>\gamma})$, then $$\tilde A_\gamma=\bigcup_{\alpha\in\Gamma_{\ge\gamma}}A_\alpha=\bigcup_{\alpha\in\Gamma_{>\gamma}}A_\alpha$$
because of the $\inf$-continuity of the system $\A$. Given any point $a\in \tilde A_\gamma$, find $\alpha\in\Gamma_{>\gamma}$ such that $a\in A_\alpha$.
Since $B\not\ni \gamma=\inf B$, there is a point $\beta\in B\cap(\gamma,\alpha)$. Now the definition of $\tilde A_\beta$ implies that $a\in A_\alpha\subset\tilde A_\beta\subset \bigcup_{\delta\in B}\tilde A_\delta$.
\smallskip

1b) If $\Gamma\ni\gamma<\inf(\Gamma_{>\gamma})$, then we can find $\beta\in B\cap (\gamma,\inf(\Gamma_{>\gamma}))$ and conclude that $\tilde A_\gamma=A_\gamma=\tilde A_\beta\subset\bigcup_{\delta\in B}\tilde A_\delta$.
\smallskip

1c) If $\Gamma\not\ni\gamma<\inf(\Gamma_{>\gamma})$, then $\sup(\Gamma_{<\gamma})\notin\Gamma$. Choose any point $\beta\in B\cap(\gamma,\inf(\Gamma_{>\gamma}))$ and observe that $\Gamma_{\ge\beta}=\Gamma_{\ge\gamma}$, $\sup(\Gamma_{<\beta})=\sup(\Gamma_{<\gamma})\notin\Gamma$ and thus $$\tilde A_\gamma=\bigcup_{\alpha\in\Gamma_{\ge\gamma}}A_\alpha=\bigcup_{\alpha\in\Gamma_{\ge\beta}}A_\alpha=\tilde A_\beta\subset\bigcup_{\delta\in B}\tilde A_\delta.$$
\smallskip

2. $\alpha=\sup(\Gamma_{<\gamma})\in\Gamma$ and $\gamma<\inf(\Gamma_{\ge\gamma})$, in which case $\tilde A_\gamma=\A_\alpha$. Since $\inf B=\gamma\notin \Gamma$, there is a point $\beta\in B\cap (\gamma,\inf(\Gamma_{\ge\gamma})$.
\smallskip

2a) If $\gamma\in\Gamma$, then $\sup(\Gamma_{<\beta})=\gamma\in\Gamma$ and thus
$\tilde A_\gamma=A_\gamma=\tilde A_\beta\subset\bigcup_{\delta\in B}A_\delta.$
\smallskip

2b) If $\gamma\notin\Gamma$, then $\Gamma_{<\beta}=\Gamma_{<\gamma}$ and thus $\tilde A_\gamma=A_\alpha=\tilde A_\beta\subset\bigcup_{\delta\in B}A_\delta.$
\end{proof}

In the following theorem for every subset $\Gamma\subset\IR$ we introduce a model $\Class_\Gamma$-absorbing system $\Sigma_\Gamma$ in the Hilbert cube $Q^\IQ$.

\begin{theorem}\label{model} For every $\Gamma\subset\IR$ the $\Gamma$-system $\Sigma_\Gamma=\la Q^\IQ,Q^{\IQ_{\le\gamma}}\times B(Q^{\IQ_{>\gamma}})\ra_{\gamma\in\Gamma}$  is $\Class_\Gamma$-absorbing and hence is homeomorphic to any other $\Class_\Gamma$-absorbing $\Gamma$-system $\la X,X_\gamma\ra_{\gamma\in\Gamma}$ in a Hilbert cube $X$.
\end{theorem}

\begin{proof} First we check that the system $\Sigma_\Gamma$ is strongly $\Class_\Gamma$-universal.

We start defining a suitable metric on the Hilbert cube $Q^\IQ$.
Let $\nu:\IQ\to(0,1)$ be any vanishing function, which means that for every $\e>0$ the set $\{q\in \IQ:\nu(q)\ge\e\}$ is finite. Take any metric $d$ generating the topology
of the Hilbert cube $Q$ and consider the metric
$$\rho((x_q),(y_q))=\max_{q\in \IQ} \nu(q)\cdot d(x_q,y_q)$$ on the Hilbert cube $Q^\IQ$.

In order to prove the strong $\Class_\Gamma$-universality of the system $\Sigma_\Gamma$,  fix a $\Gamma$-system $\mc A=\la A,A_\gamma\ra_{\gamma\in\Gamma}\in\Class_\Gamma$ and a map $f:A\to Q^\IQ$ that restricts to a $Z$-embedding of some closed subset $K\subset A$ such that $(f|K)^{-1}(\Sigma_\Gamma)=K\cap \mc A$. Given $\e>0$, we need to construct a  $Z$-embedding $\tilde f:A\to Q^\IQ$ such that $\rho(\tilde f,f)<\e$, $\tilde f|K=f|K$ and $\tilde f^{-1}(\Sigma_\gamma)=A_\gamma$ for all $\gamma\in\Gamma$.

By Lemma~\ref{enlarge}, the $\Gamma$-system $\A$ extends to an $\IR$-system $\tilde\A=\la A,A_\gamma\ra_{\gamma\in\IR}\in\Class_\IR$. We shall construct a  $Z$-embedding $\tilde f:A\to Q^\IQ$ such that $\rho(\tilde f,f)<\e$, $\tilde f|K=f|K$ and $\tilde f^{-1}(\Sigma_\gamma)\setminus K=A_\gamma\setminus K$ for all $\gamma\in\IR$.

For every $q\in \IQ$ let  $\pr_q:Q^\IQ\to Q$ denote the coordinate projection.
Since $f(K)$ is a $Z$-set in $Q^\IQ$, we can approximate the map $f$ by a map $f':A\to Q^\IQ$ such that $\rho(f',f)<\e/2$, $f'|K=f|K$ and $f'(A\setminus K)\cap f'(K)=\emptyset$.
Using the strong $\Class_{\{0\}}$-universality of the pair $(Q,B(Q))$,
for each $q\in \IQ$ we can approximate the map $\pr_q\circ f':A\to Q$ by a map $\tilde f_q:A\to Q$ such that
\begin{itemize}
\item[a)] $d(\tilde f_q(x),\pr_q\circ f'(x))\le \frac\e2 \rho(f'(x),f(K))$ for all $x\in A$;
\item[b)] $\tilde f_q|A\setminus K$ is injective;
\item[c)] $\tilde f_q(A\setminus K)$ is a $\sigma Z$-set in $Q$;
\item[d)] $\tilde f_q^{-1}(B(Q))\setminus K=A_q\setminus K$.
\end{itemize}

Now consider the diagonal product $\tilde f=(\tilde f_q)_{q\in \IQ}:A\to Q^\IQ$ of the maps $\tilde f_q$, $q\in \IQ$. It follows from (a) that $\tilde f|K=f'|K=f|K$, $\rho(\tilde f,f)\le \rho(\tilde f,f')+\rho(f',f)<\e$ and $\tilde f(A\setminus K)\cap f(K)=\emptyset$.
Combining this fact with (b) we conclude that the map $\tilde f:A\to Q^\IQ$ is injective and hence an embedding. It follows from (c) that $\tilde f(A)$ is a $\sigma Z$-set in $Q^\IQ$ and hence a $Z$-set, see \cite[6.2.2]{jan}. Therefore, $\tilde f$ is a $Z$-embedding approximating the map $f$.
\medskip

It remains to check that $\tilde f^{-1}(\Sigma_\gamma)=A_\gamma$ for every $\gamma\in\Gamma$. Since $\tilde f^{-1}(\Sigma_\gamma)\cap K=(f|K)^{-1}(\Sigma_\gamma)=K\cap A_\gamma$, it suffices to check that
$\tilde f^{-1}(\Sigma_\gamma)\setminus K=A_\gamma\setminus K$.

It follows that
$$
\begin{aligned}
\tilde f^{-1}(\Sigma_\gamma)\setminus K&=\tilde f^{-1}(Q^{\IQ_{\le\gamma}}\times B(Q^{\IQ_{>\gamma}}))\setminus K=\\&=\bigcup_{q\in \IQ_{>\gamma}}\tilde f_q^{-1}(B(Q))\setminus K
=\bigcup_{q\in \IQ_{>\gamma}}A_q\setminus K=A_\gamma\setminus K.
\end{aligned}
$$
The last equality follows from the $\inf$-continuity of the $\IR$-system $\widetilde{\mc A}=\la A,A_\gamma\ra_{\gamma\in\IR}$ because $\gamma=\inf \IQ_{>\gamma}$.
This completes the proof of the strong $\Class_\Gamma$-universality of the system $\Sigma_\Gamma$.
\smallskip

It remains to check that the $\Gamma$-system $\Sigma_\Gamma$ satisfies the second condition of Definition~\ref{absorb} of a $\Class_\Gamma$-absorbing system. It is clear the $\Gamma$-system $\Sigma_\Gamma$ is $\sigma$-compact and decreasing.   To show that it is $\inf$-continuous, take any subset $B\subset \Gamma$ with $\gamma=\inf B\in\Gamma$. If $\gamma\in B$, then $\Sigma_\gamma\supset\bigcup_{\beta\in B}\Sigma_\beta\supset \Sigma_\gamma$. So, we assume that $\gamma\notin B$.  Since the $\Gamma$-system $\la Q^\IQ, \Sigma_\gamma\ra_{\gamma\in\Gamma}$ is decreasing, we get $\Sigma_\gamma\supset\bigcup_{\beta\in B}\Sigma_\beta$. To prove the reverse inclusion, take any point
$(x_q)_{q\in \IQ}\in\Sigma_\gamma=Q^{\IQ_{\le\gamma}}\times B(Q^{\IQ_{>\gamma}})$ and observe that $x_q\in B(Q)$ for some $q\in\IQ_{>\gamma}$. Since $\gamma=\inf B$ the half-interval $[\gamma,q)$ contains a point $\beta\in B$.

Then $(x_q)_{q\in \IQ}\in Q^{\IQ_{\le\beta}}\times B(Q^{\IQ_{>\beta}})$ and thus $(x_q)_{q\in \IQ}\in \Sigma_\beta\subset\bigcup_{\alpha\in B}\Sigma_\alpha$.
Therefore, $\Sigma_\Gamma\in\Class_\Gamma$.
\smallskip

Since each space $\Sigma_\gamma$, $\gamma\in\IQ$, is a $\sigma Z$-set in $Q^\IQ$, so is the countable union $\bigcup_{\gamma\in \IQ}\Sigma_\gamma=\bigcup_{\gamma\in\IR}\Sigma_\gamma$. So, we can find a sequence $\la Z_n\ra_{n\in\w}$ of $Z$-sets in $Q^\IQ$ such that $$\bigcup_{n\in\w}Z_n=\bigcup_{\gamma\in \IQ}\Sigma_\gamma.$$ It follows from $\Sigma_\Gamma\in\Class_\Gamma$ that $Z_n\cap\Sigma_\Gamma\in\Class_\Gamma,$
which completes the proof of the $\Class_\Gamma$-absorbing property of the system $\Sigma_\Gamma$.
\smallskip

By the Uniqueness Theorem~\ref{unique}, each $\Class_\Gamma$-absorbing system $\la
X,X_\gamma\ra_{\gamma\in\Gamma}$ in a Hilbert cube $X$ is homeomorphic to the $\Class_\Gamma$-absorbing
$\Gamma$-system $\Sigma_\Gamma$.
\end{proof}

\section{Strongly universal systems of hyperspaces}\label{s2}

In this section we establish an important Theorem~\ref{cauty} detecting strongly ${\mathfrak C}_\Gamma$-universal $\Gamma$-systems in hyperspaces. In this section, $\Gamma$ is any set and $\mathfrak C_\Gamma$ is a class of $\Gamma$-systems.

By the {\em hyperspace} of a topological space $X$ we understand the space $2^X$ of nonempty compact subsets of $X$ endowed with the Vietoris topology. This topology is generated by the sub-base consisting of the sets
$$
\la V\ra=\{K\in 2^X:K\subset V\}\mbox{ and }\la X,V\ra=\{K\in 2^X:K\cap V\ne\emptyset\}$$
where $V$ is an open subset of $X$. If the topology of $X$ is generated by a metric $d$, then the Vietoris topology on $2^X$ is generated by the Hausdorff metric
$d_H(A,B) = \max\{\max_{a\in A}d(a,B),\max_{b\in B}d(b,A)\}$.

In the sequel by $2_{<\w}^X$ we shall denote the subspace of $2^X$ consisting of finite non-empty subsets of $X$. By \cite{glad}, \cite[8.4.3]{jan} for a non-degenerate Peano continuum $X$ the subset $2_{<\w}^X$ is  homotopy dense in $2^X$. 

We recall that a subset $A$ of a topological space $X$ is {\em homotopy dense} if there is a homotopy $h:X\times [0,1]\to X$ such that $h(x,0)=x$ and $h(x,t)\in  A$ for all $x\in X$ and $t\in(0,1]$.
\smallskip

We define a subspace $\mathcal H\subset 2^X$ to be {\em finitely additive} if 
\begin{itemize}
\item $A\cup F\in\mathcal H$ for any $A\in\mathcal H$ and any finite subset $F\subset X$; 
\item $A\sqcup B\in\mathcal H$ for any disjoint sets $A,B\in \mathcal H$. 
\end{itemize}

The first condition implies that each finite subset of $X$ belongs to the family
$$\add(\mathcal H)=\{A\in 2^X:\forall B\in \mathcal H\;\; A\cup B\in \mathcal H\}.$$

For a $\Gamma$-system $\mathscr H=\la 2^X,\HH_\gamma\ra_{\gamma\in\Gamma}$ the intersection
$$\add(\mathscr H)=\bigcap_{\gamma\in\Gamma}\add(\HH_\gamma)\cap \add(2^X\setminus\HH_\gamma)$$ will be called {\em the additive kernel} of $\mathscr H$.

For example, the additive kernel of the $\w$-system $\la 2^X,\dim_{\le n}(X)\ra_{n\in\w}$ is equal to the subspace $\dim_{\le 0}(X)$ of all zero-dimensional compact subsets of $X$. The additive kernel of the $[0,\infty)$-system
$\la 2^X,\HD_{\le\gamma}(X)\ra_{\gamma\in[0,\infty)}$ is equal to the subspace $\HD_{\le 0}(X)\subset 2^X$ consisting of subsets of $X$ with Hausdorff dimension zero.

The following technical theorem was implicitly proved by R.Cauty in \cite{coty}.

\begin{theorem}\label{cauty} Let $X$ be a non-degenerate Peano continuum. A $\Gamma$-system $\mathscr H=\la 2^X,\mathcal {\mathcal H_\gamma}\ra_{\gamma\in\Gamma}$ is strongly ${\mathfrak C}_\Gamma$-universal if:
\begin{enumerate}
\item[1)] for every $\gamma\in\Gamma$ the subspaces $\HH_\gamma$ and $2^X\setminus\HH_\gamma$ are finitely additive;
\item[2)] for every non-empty open set $U\subset X$ there is a map $\xi:Q\to 2^U\cap \add(\mathscr H)$ such that for any distinct points $x,x'\in Q$ the symmetric difference $\xi(x)\triangle\xi(x')$ is infinite;
\item[3)] for any non-empty open set $U\subset X$ and any $\Gamma$-system $\mc C=\la C,C_\gamma\ra_{\gamma\in\Gamma}\in{\mathfrak C}_\Gamma$ there is a map $\varphi:C\to 2^U$ such that $\varphi^{-1}(\mathscr H)=\mathscr C$.
\end{enumerate}
\end{theorem}

\section{The strong $\Class_\Gamma$-universality of $\Gamma$-systems of hyperspaces}

In this section,  we detect
strongly $\Class_\Gamma$-universal systems of the form $\la 2^X,\D_{>\gamma}(X)\ra_{\gamma\in\Gamma}$
where $\Gamma\subset[0,\infty)$ and $\D:2^X_*\to[0,\infty]$ is a dimension function defined on the hyperspace of a non-degenerated Peano continuum $X$.

First we establish one property of  dimension functions which is formally stronger that the $\w$-additivity.

\begin{lemma}\label{omega} Let $X$ be a metrizable compact space without isolated points and $\D:2_*^X\to[0,\infty]$ be a dimension function. For every non-empty open set $U\subset X$ there is a disjoint sequence $\la U_n\ra_{n\in\w}$ of non-empty open sets of $U$ such that
\begin{enumerate}
\item $\la U_n\ra_{n\in\w}$ converges to some point $x_\infty\in U$, which means that each neighborhood $O(x_\infty)$ contains all but finitely many sets $U_n$;
\item for any compact subsets $K_n\subset U_n$, $n\in\w$, the set $K_\infty=\{x_\infty\}\cup\bigcup_{n\in\w}K_n$ is compact and has dimension $\D(K_\infty)\le\sup_{n\in\w}\D(K_n)$.
\end{enumerate}
\end{lemma}

\begin{proof} Take any non-empty open subset $V\subset X$ with $\cl(V)\subset U$. The $\w$-additivity of the dimension function $\D$ yields a sequence $\la V_n\ra_{n\in\w}$ of open subsets of $V$ such that for any compact subset $K\subset\cl(\bigcup_{n\in\w}V_n)$ has dimension  $\D(K)\le\sup_{n\in\w}\D(K\cap\overline{V}_n)$.

Replacing the sets $V_n$ by their suitable subsets, we can assume that $\diam(V_n)\to 0$ as $n\to\infty$. In
each set $V_n$ pick a point $x_n$. Since the space $X$ has no isolated point, we can choose the points $x_n$,
$n\in\w$, to be pairwise distinct. Next, replacing the sets $V_n$ by small neighborhoods of the points $x_n$, we
can make the sets $V_n$, $n\in\w$, pairwise disjoint. By the compactness of $X$, the sequence $\la
x_n\ra_{n\in\w}$ contains a subsequence $\la x_{n_k}\ra_{k\in\w}$ that converges to some point $x_\infty\in
\cl(V)\subset U$. Since $\diam(V_{n_k})\to 0$, the sequence $\la V_{n_k}\ra_{k\in\w}$ also converges to
$x_\infty$.

It is clear that the sets $U_k=V_{n_k}$, $k\in\w$, have the desired properties.
\end{proof}

Now we are able to prove the principal ingredient in the proof of Theorem~\ref{main}.
Below $\Gamma\subset[0,\infty)$ and $\Class_\Gamma$ stands for the class of $\inf$-continuous $\sigma$-compact $\Gamma$-systems.

\begin{theorem}\label{suh}  Let $X$ be a non-degenerate Peano continuum, $\D:2_*^X\to[0,\infty]$ be a dimension function, and $\Gamma\subset[0,\infty)$. The $\Gamma$-system $\la 2^X,\D_{>\gamma}(X)\ra_{\gamma\in\Gamma}$ is strongly $\Class_\Gamma$-universal if and only if
each non-empty open set $U\subset X$ for every $\gamma\in\Gamma$ contains a compact subset $K\subset U$ with $\D(K)\in (\gamma]_\Gamma$.
\end{theorem}

\begin{proof} To prove the ``only if'' part, assume that the system $\mc D=\la 2^X,\D_{>\gamma}(X)\ra_{\gamma\in\Gamma}$ is strongly $\Class_\Gamma$-universal.

Fix any non-empty open set $U\subset X$ and an element $\gamma\in\Gamma$. We need to find a compact subset $K\subset U$ with $\D(K)\in(\gamma]_\Gamma$.

Let $A=\{a\}$ be any singleton and put $A_\alpha=A$ for all $\alpha<\gamma$ and $A_\alpha=\emptyset$ for all $\alpha>\gamma$. Put also $A_\gamma=\emptyset$ if $\gamma=\inf(\Gamma_{>\gamma})$ and $A_\gamma=A$ otherwise.

Observe that the so-defined $\Gamma$-system $\mc A=\la A,A_\gamma\ra_{\gamma\in\Gamma}$ belongs to the class $\Class_\Gamma$. Now using the strong $\Class_\Gamma$-universality of the $\Gamma$-system $\mc D$, find a map $f:A\to 2^U$ such that $f^{-1}(\mc D)=\mc A$.

We claim that the compact subset $K=f(a)\subset U$ has dimension $\D(K)\in(\gamma]_\Gamma$.
To prove this inclusion, consider the three cases from the definition of the set $(\gamma]_\Gamma$.
\smallskip

(i) If $\gamma<\inf(\Gamma_{>\gamma})$, then $a\in A_\gamma$ and hence $K=f(a)\in \D_{>\gamma}(X)$ and $\gamma<\D(K)$.
On the other hand, for every $\alpha\in\Gamma_{>\gamma}$ we get $a\notin A_\alpha=\emptyset$ and thus $K=f(a)\in 2^X\setminus D_{>\alpha}(X)=D_{\le \alpha}(X)$ and  $\D(K)\le\alpha$, which implies $\D(K)\le\inf(\Gamma_{>\gamma})$.
Consequently, $\D(K)\in(\gamma,\inf(\Gamma_{>\gamma})]=(\gamma]_\Gamma$.

(ii) $\Gamma\ni\sup(\Gamma_{<\Gamma})<\gamma=\inf(\Gamma_{>\gamma})$. In this case $a\notin A_\gamma=\emptyset$ and thus $K=f(a)\in \D_{\le\gamma}(X)$. On the other hand, $a\in A_\alpha$ where $\alpha=\sup(\Gamma_{<\gamma})<\gamma$ and hence $K=f(a)\in\D_{>\alpha}(X)$. Conseqeuntly, $\D(K)\in (\sup(\Gamma_{<\gamma}),\gamma]=(\gamma]_\Gamma$.

(iii) If $\gamma=\inf(\Gamma_{>\gamma})$ and $\sup(\Gamma_{<\gamma})$  is equal $\gamma$ or does not belongs to
$\Gamma$, then for every $\alpha\in \Gamma_{<\gamma}$, we get $a\in A_\alpha$ and thus
$K=f(a)\in\D_{>\alpha}(X)$ ad $\D(K)>\alpha$. Consequently, $\D(K)\ge\sup(\Gamma_{<\gamma})$. On the other hand,
$a\notin A_\gamma=\emptyset$ implies $K=f(a)\in\D_{\le\gamma}(X)$ and thus
$\D(K)\in[\sup(\Gamma_{<\gamma}),\gamma]=(\gamma]_\Gamma$.
\medskip

To prove the ``only if'' part, assume that for every non-empty open set $U\subset X$ and every $\gamma\in\Gamma$ there is a compact subset $K\subset U$ with $\D(K)\in(\gamma]_\Gamma$.

The strong $\Class_\Gamma$-universality of the system  $\mc D$ will follow as soon as we check the conditions (1)--(3) of Theorem~\ref{cauty} for the class $\Class_\Gamma$.
\smallskip

1. The monotonicity of the dimension function $\D$ implies that the subspace $\D_{>\gamma}(X)$ of $2^X$ is finitely additive. The finite additivity of the complement $\D_{\le\gamma}(X)=2^X\setminus\D_{>\gamma}(X)$ follows from the finite additivity of the dimension function $\D$.
\smallskip

2. To establish the condition (2) of Theorem~\ref{cauty}, fix any non-empty open set $U\subset X$. Lemma~\ref{omega} yields a sequence  $\la U_n\ra_{n\in\w}$ of non-empty open subsets of $U$ that converge to some point $x_\infty\in U$ and has the property that for any compact subsets $K_n\subset U_n$ the set $K=\{x_\infty\}\cup\bigcup_{n\in\w}K_n$ is compact and has dimension $\D(K)\le\sup_{n\in\w}\D(K_n)$. Each set $U_n$ contains a topological copy of the interval $[0,1]$, so we can find a topological embedding $\xi_n:[-1,1]\to U_n$.

Let $\nu:\w\to\w$ be any function such that the preimage $\nu^{-1}(n)$ of every $n\in\w$ is infinite. Define a map $\xi:Q\to 2^U$ assigning to each $\vec t=\la t_n\ra_{n\in\w}\in Q$ the compact subset
$$\xi(\vec t)=\{x_\infty\}\cup\{\alpha_n(t_{\nu(n)}):n\in\w\}$$ of $U$ having a unique non-isolated point $x_\infty$. The equality $\D(\emptyset)=0$ and the finite additivity of the dimension function $\D$ implies that $\D(F)=0$ for each finite subset $F\subset X$. The choice of the sequence $\la U_n\ra$ guarantees that $\D(\xi(\vec t))=0$ and thus $$\xi(Q)\subset \D_{\le 0}(X)\subset\add(\mc D).$$ The choice of the function $\nu$ guarantees that $\xi(\vec t)\triangle\xi(\vec u)$ is infinite for any distinct vectors $\vec t,\vec u\in Q$.
\smallskip

3. To check the condition (3) of Theorem~\ref{cauty}, fix any non-empty open set $U\subset X$ and a $\Gamma$-system $\mc A=\la A,A_\gamma\ra_{\gamma\in\Gamma}\in\Class_\Gamma$.
Each set $A_\gamma$, $\gamma\in\Gamma$, being $\sigma$-compact, can be written as the countable union $A_\gamma=\bigcup_{n\in\w}A_{\gamma,n}$ of an increasing sequence $\la A_{\gamma,n}\ra_{n\in\w}$ of compact subsets of $A$.
Let $D$ be a countable subset of $\Gamma$ meeting each half-interval $[\gamma,\gamma+\e)$ where $\gamma\in\Gamma$ and $\e>0$.

Apply Lemma~\ref{omega} to find a disjoint family $\la U_d\ra_{d\in D}$ of non-empty open subsets of $U$ such that
\begin{itemize}
\item $\la U_d\ra_{d\in D}$ converges to some point $x_\infty\in U$ in the sense that each neighborhood $O(x_\infty)$ contains all but finitely many sets $U_d$, $d\in D$;
\item for any compact sets $K_d\subset U_d$ the set $K=\{x_\infty\}\cup\bigcup_{d\in D}K_d$ is compact and has dimension $\D(K)\le\sup_{d\in D}\D(K_d)$.
\end{itemize}

For every $d\in D$ use Lemma~\ref{omega} once more and find a disjoint family $\la U_{d,n}\ra_{n\in\w}$ of non-empty open subsets of $U$ such that
\begin{itemize}
\item $\la U_{d,n}\ra_{n\in\w}$ converges to some point $x_d\in U_d$;
\item for any compact sets $K_{n}\subset U_{d,n}$ the set $K_d=\{x_d\}\cup\bigcup_{n\in\w}K_{n}$ is compact and has dimension $\D(K_d)=\sup_{n\in\w}\D(K_n)$.
\end{itemize}

By our assumption, for every $d\in D$ and $n\in\w$ we can find a compact subset
$K_{d,n}\subset {U_{d,n}}$ with $\D(K_{d,n})\in(d]_\Gamma$.
Using the homotopical density of the subspace $2_{<\w}^X $ of finite subsets in $2^X$, construct a map
$\kappa_{d,n}:A\to 2^X$ such that $\kappa_{d,n}(a)=K_{d,n}$ for every $a\in A_{d,n}$ and $\kappa_{d,n}(a)$ is a finite subset of $U_{d,n}$ for every $a\in A\setminus A_{d,n}$.

Now for every $a\in A$ and $d\in D$ consider the compact subset $$\kappa_{d}(a)=\{x_d\}\cup\bigcup_{n\in\w}\kappa_{d,n}(a)\subset U_d$$
having dimension
$$\D(\kappa_d(a))=\sup_{n\in\w}\D(\kappa_{d,n}(a)).$$

The choice of the sequence $\la U_d\ra_{d\in D}$ ensures that
$$\kappa(a)=\{x_\infty\}\cup\bigcup_{d\in D}\kappa_d(a)$$
is a compact subset of $U$ with dimension
$$\D(\kappa(a))=\sup_{d\in D}\D(\kappa_d(a))=\sup\{\D(\kappa_{d,n}(a)):d\in D,\;n\in\w\}.$$

It is easy to prove that the map
$$\kappa:A\to 2^U,\;\; \kappa:a\mapsto \kappa(a),$$
is continuous. It remains to check that $\kappa^{-1}(\D_{>\gamma}(X))=A_\gamma$ for all $\gamma\in \Gamma$.
\smallskip

If $a\in A\setminus A_\gamma$, then for every $d\ge \gamma$ in $D$ the inclusion $a\in A\setminus A_d$ implies  $\kappa_{d,n}(a)\in 2_{<\w}^X $. In this case $\D(\kappa_d(a))\le \sup_{n\in\w}\kappa_{d,n}(a)=0\le\gamma$.
On the other hand, for every $d<\gamma$ the inclusions $D(K_{d,n})\in(d]_\Gamma\subset[0,\gamma]$, $n\in\w$, and the choice of the sequence $\la U_{d,n}\ra_{n\in\w}$ imply  $\D(\kappa_d(a))\le \sup_{n\in\w}D(\kappa_{d,n}(a))\le\gamma$.

Now the choice of the sequence $\la U_d\ra_{d\in D}$ guarantees that $$\D(\kappa(a))\le \sup_{d\in D}\D(\kappa_d(a))\le \gamma$$ and hence
$$\kappa(a)\in\D_{\le \gamma}(X)=2^X\setminus\D_{>\gamma}(X).$$
\smallskip

Now assume that $a\in A_\gamma$ and hence $a\in A_{\gamma,n}$ for some $n\in\w$. If $\gamma<\inf(\Gamma_{>\gamma})$, then $\gamma\in D$ and $\D(K_{\gamma,n})\in (\gamma]_\Gamma=(\gamma,\inf(\Gamma_{>\gamma})]$. Since $K_{\gamma,n}\subset \kappa(a)$, we conclude that $\D(\kappa(a))\ge\D(K_{\gamma,n})>\gamma$ and thus $\kappa(a)\in\D_{>\gamma}(X)$.

Next, assume that $\gamma=\inf(\Gamma_{>\gamma})$. In this case $\gamma=\inf(D_{>\gamma})$ and hence $A_\gamma=\bigcup_{d\in D_{>\gamma}}A_d$. It follows that $a\in A_{d,n}$ for some $d\in D_{>\gamma}$ and $n\in\w$. Since $\kappa(a)\supset K_{d,n}$ and $\D(K_{d,n})\in(d]_\Gamma\subset (\gamma,+\infty)$, we conclude that $\D(\kappa(a))\ge\D(K_{d,n})>\gamma$. So, again $\kappa(a)\in\D_{>\gamma}(X)$.
\end{proof}

The following characterization theorem implies Theorem~\ref{main} announced in the Introduction.

\begin{theorem}\label{absh} Let $X$ be a topological space, $\D:2_*^X\to[0,\infty]$ be a dimension function, and $\Gamma\subset[0,\infty)$ be a subset. The $\Gamma$-system $\la 2^X,\D_{>\gamma}(X)\ra_{\gamma\in\Gamma}$ is homeomorphic to the model $\Class_\Gamma$-absorbing $\Gamma$-system $\la Q^\IQ,Q^{\IQ_{\le\gamma}}\times B(Q^{\IQ_{>\gamma}})\ra_{\gamma\in\Gamma}$ if and only if
\begin{enumerate}
\item $X$ is a non-degenerate Peano continuum,
\item each space $\D_{>\gamma}(X)$, $\gamma\in\Gamma$, is $\sigma$-compact, and
\item each non-empty open set $U\subset X$ for every $\gamma\in\Gamma$ contains a compact subset $K\subset U$ with $D(K)\in(\gamma]_\Gamma$.
\end{enumerate}
\end{theorem}

\begin{proof} To prove the ``only if'' part, assume that the $\Gamma$-system $\mathscr D=\la 2^X,\D_{>\gamma}(X)\ra_{\gamma\in\Gamma}$ is homeomorphic to the model $\Gamma$-system $\Sigma_\Gamma=\la Q^\IQ,Q^{\IQ_{\le\gamma}}\times B(Q^{\IQ_{>\gamma}})\ra_{\gamma\in\Gamma}$. Since $2^X$ is homeomorphic to $Q^\IQ$, we may apply the Curtis-Shori Theorem \cite{CS} and conclude that $X$ is a non-degenerate Peano continuum.

Since each space $\Sigma_\gamma=Q^{\IQ_{\le\gamma}}\times B(Q^{\IQ_{>\gamma}})$, $\gamma\in\Gamma$, is $\sigma$-compact, so is its topological copy $\D_{>\gamma}(X)$.

The $\Gamma$-system $\mathscr D$, being homeomorphic to the model $\Class_\Gamma$-absorbing $\Gamma$-system $\Sigma_\Gamma$, is strongly $\Gamma$-universal. Now Theorem~\ref{suh} guarantees that for every $\gamma\in\Gamma$ each non-empty open subset $U\subset X$ contains a compact subset $K\subset U$ with $\D(K)\in(\gamma]_\Gamma$.
\smallskip

Next, we prove the ``if'' part. Assume that the conditions (1)--(3) are satisfied. We shall prove that the $\Gamma$-system $\mathscr D$ is $\Class_\Gamma$-absorbing. By the Curtis-Schori Theorem \cite{CS}, the hyperspace $2^X$ is homeomorphic to the Hilbert cube $Q$. By Theorem~\ref{suh}, the $\Gamma$-system $\mathscr D$ is strongly $\Class_\Gamma$-universal. It is clear that this $\Gamma$-system is $\inf$-continuous. By the condition (2), it is $\sigma$-compact.  Hence $\mathscr D\in\Class_\Gamma$.

Let $D\subset\Gamma$ be countable subset that meets each half-interval $[\gamma,\gamma+\e)$ where
$\gamma\in\Gamma$ and $\e>0$. It follows that $\bigcup_{\gamma\in\Gamma}\D_{>\gamma}(X)=\bigcup_{\gamma\in
D}\D_{>\gamma}(X)\subset\D_{>0}(X)$ is a $\sigma Z$-set in $2^X$, being  a $\sigma$-compact subset of $2^X$ that
has empty intersection with the homotopy dense subset $2^X_{<\w}\subset\D_{\le0}(X)$ on $2^X$. So, we can find a
countable sequence $\la Z_n\ra_{n\in\w}$ of $Z$-sets in $2^X$ such that $\bigcup_{n\in\w}Z_n\supset
\bigcup_{\gamma\in\Gamma}\D_{>\gamma}(X)$. Since $\mathscr D\in\Class_\Gamma$, we get $Z_n\cap\mathscr
D\in\Class_\Gamma$ for all $n\in\w$. This completes the proof of the $\Class_\Gamma$-absorbing property of the
$\Gamma$-system $\mathscr D$. Since $2^X$ is homeomorphic to the Hilbert cube, Theorem~\ref{model} ensures that
$\mathscr D$ is homeomorphic to the model $\Gamma$-system $\Sigma_\Gamma$.
\end{proof}

\section{Mean Value Theorem for Hausdorff dimension}\label{s6}

In this section we shall prove Theorem~\ref{dimH}. 
First, we recall shortly the definitions of the Hausdorff measure and dimension. Given a complete separable metric
space $E$ and two non-negative real numbers $s$, $\varepsilon$, consider the number
$$
{{\mathcal H}}_\varepsilon^s(E) = \inf\limits_{{\mathcal B}}\sum\limits_{B\in{{\mathcal B}}} ({\rm d}{\rm i}{\rm
a}{\rm m} B)^s,
$$
where infimum is taken over all $\e$-covers ${\mathcal B}$ of $E$, i.e. cover of $E$ by sets of diameter $\le\e$.
Since $X$ is separable, we can restrict ourselves by countable covers by closed subsets of diameter $\le\e$.

The limit ${{\mathcal H}}^s(E)=\lim\limits_{\varepsilon \rightarrow 0}{{\mathcal H}}_ \varepsilon^s(E)$ is
called the $s$-dimensional Hausdorff measure of $E$. It is known that there is a unique finite or infinite
number $\dim_H(E)$ called the \textit{Hausdorff dimension} of $E$ and denoted by $\dim_H(E)$ such that
${{\mathcal H}}^s(E)=\infty$ for all $s < \dim_H(E)$ and ${{\mathcal H}}^s(E)=0$ for all $s>\dim_H(E)$, see
\cite{edg}, \cite{fal}.

Let $(X,d)$ be a separable complete metric space. Theorem~\ref{dimH} will be proved as soon as for every positive  real number $s<\dim_H(X)$ we shall find a compact subset $K\subset X$ with Hausdorff dimension $\dim_H(K)=s$.

It follows from $s<\dim_H(E)$ that ${\mathcal H}^s(E)=\infty$ and there exists $0<\delta<1$ with
$$
{\mathcal H}_\delta^s(E)=k_0>\frac{\delta^s}{2^{s-1}}=\left(\frac\delta2\right)^s+\left(\frac\delta4\right)^s+
\left(\frac\delta8\right)^s+\cdots.\eqno(0)
$$

We define inductively a decreasing sequence $\{E_i\}_{i=1}^\infty$ of closed subsets of $E$. Let $E_1=E$.
Consider ${\mathcal H}_{\delta/2}^s(E_1)$. Two cases are possible (taking into account the definition of
Hausdorff measure):
\begin{itemize}
    \item ${\mathcal H}_{\delta/2}^s(E_1)=k_0$. In this case we take $E_2=E_1$.
    \item ${\mathcal H}_{\delta/2}^s(E_1)>k_0$. Therefore we can choose a closed $\delta/2$-cover $\{U_1,\ldots,U_{m_1},\ldots\}$ of the set $E_1$,
     (without loss of generality assume that this cover is ordered so that $\diam(U_{i+1})\le\diam(U_i)$ for all $i$), such that
$$
{\mathcal H}_{\delta/2}^s(E_1)\le\sum_{i}(\diam(U_i))^s<{\mathcal H}_{\delta/2}^s(E_1)+(\delta/2)^s.\eqno(1)
$$
Find a finite number $m_1$ such that
$$
k_0\le\sum_{i=1}^{m_1}(\diam(U_i))^s\le k_0+(\delta/2)^s. \eqno(2)
$$
Then take $E_2=\bigcup\limits_{i=1}^{m_1}E_1\cap U_i$.
\end{itemize}

Now we need to estimate ${\mathcal H}_{\delta/2}^s(E_2)$(obviously the second case is interesting). For this we
put $E_2'=\bigcup\limits_{i>m_1}E_1\cap U_i$ and note that
$$
{\mathcal H}_{\delta/2}^s(E_1)\le{\mathcal H}_{\delta/2}^s(E_2)+{\mathcal H}_{\delta/2}^s(E_2'). \eqno(3)
$$
On the other hand
$$
\sum_{i}(\diam(U_i))^s=\sum_{i\le m_1}(\diam(U_i))^s + \sum_{i>m_1}(\diam(U_i))^s. \eqno(4)
$$

Consider the real numbers  $\e_1=\sum_{i}(\diam(U_i))^s-{\mathcal H}_{\delta/2}^s(E_1)$,
$\e_2=\sum_{i=1}^{m_1}(\diam(U_i))^s-{\mathcal H}_{\delta/2}^s(E_2)$,
$\e_2'=\sum_{i>m_1}(\diam(U_i))^s-{\mathcal H}_{\delta/2}^s(E_2')$ and observe that 
$0\le\varepsilon_1<(\delta/2)^s$ by (1), and $\varepsilon_2,\varepsilon_2'\ge 0$. Therefore  (3) and (4) yield $\varepsilon_1\ge\varepsilon_2+\varepsilon_2'$ and hence $0\le\varepsilon_2<(\delta/2)^s$. Taking into
account (2), we have:
$$
k_0-(\delta/2)^s\le{\mathcal H}_{\delta/2}^s(E_2)\le k_0+(\delta/2)^s. \eqno(5)
$$

Now denote ${\mathcal H}_{\delta/2}^s(E_2)=k_1$. From (0) and (5) it follows that $0<k_1<\infty$. By the
definition of Hausdorff measure we have $0<{\mathcal H}^s(E_2)\le\infty$, that in turn implies $\dim_H(E_2)\ge
s$. It allows us to make the following inductive step.

Consider now ${\mathcal H}_{\delta/4}^s(E_2)$. If ${\mathcal H}_{\delta/4}^s(E_2)=k_1$, then take $E_3=E_2$. If
${\mathcal H}_{\delta/4}^s(E_2)>k_1$, then similarly to the described above we find a closed $\delta/4$-cover $\{U_1,\ldots,U_{m_2},\ldots\}$ of $E_2$,
 such that
$$
{\mathcal H}_{\delta/4}^s(E_2)\le\sum_{i}(\diam(U_i))^s<{\mathcal H}_{\delta/4}^s(E_2)+(\delta/4)^s.
$$
Find a finite number $m_2$ such that
$$
k_1\le\sum_{i=1}^{m_2}(\diam(U_i))^s\le k_1+(\delta/4)^s.
$$
Let $E_3=\bigcup\limits_{i=1}^{m_2}E_2\cap U_i$. As above, we can to estimate ${\mathcal H}_{\delta/4}^s(E_3)$.
We obtain:
$$
k_1-(\delta/4)^s\le{\mathcal H}_{\delta/4}^s(E_3)\le k_1+(\delta/4)^s.
$$
Or, taking into account (5):
$$
k_0-(\delta/2)^s-(\delta/4)^s\le{\mathcal H}_{\delta/4}^s(E_3)\le k_0+(\delta/2)^s+(\delta/4)^s.
$$

Again we can state that $\dim_H(E_3)\ge s$ and continue inductive process by constructing in similar way $E_4,
E_5,\ldots E_n,\ldots$, for which we obtain in general case the estimate:
$$
k_0-(\delta/2)^s-\cdots-(\delta/2^{n-1})^s\le{\mathcal H}_{\delta/2^{n-1}}^s(E_n)\le
k_0+(\delta/2)^s+\cdots+(\delta/2^{n-1})^s. \eqno(6)
$$

It follows that $E_1\supseteq E_2\supseteq E_3\supseteq\ldots$ is a decreasing sequence of closed subsets of $X$ with compact intersection
$K=\lim\limits_{n\to\infty}E_n=\bigcap\limits_{n=1}^\infty E_n$. Using the continuity of the measure ${\mathcal H}^s$ we obtain:
$$
{\mathcal H}^s(K)=\lim_{n\to\infty}{\mathcal H}^s(E_n)=\lim_{n\to\infty}\lim_{i\to\infty}{\mathcal
H}_{\delta/2^{i-1}}^s(E_n)=\lim_{n\to\infty}{\mathcal H}_{\delta/2^{n-1}}^s(E_n).
$$
Additionally using (6) we obtain the estimate:
$$
k_0-\frac{\delta^s}{2^{s-1}}\le{\mathcal H}^s(F)\le k_0+\frac{\delta^s}{2^{s-1}}.
$$
Taking into account (0) we can state that $\dim_H(F)=s$.

\bibliographystyle{amsplain}
\bibliography{8}

\end{document}